\documentclass[12pt]{article}
\usepackage{subeqn}
\usepackage{graphicx}
\usepackage{ifpdf}
\ifpdf  \fi
\usepackage{amsfonts,amssymb}
\usepackage{mathptmx,helvet,courier,makeidx,multicol,footmisc}
\usepackage[numbers]{natbib}
\bibpunct{(}{)}{;}{a}{,}{,}
\usepackage[bookmarksnumbered=true, pdfauthor={Wen-Long Jin}]{hyperref}

\newcommand{\commentout}[1]{}

\newcommand{\ba}{\begin{array}}
        \newcommand{\ea}{\end{array}}
\newcommand{\bc}{\begin{center}}
        \newcommand{\ec}{\end{center}}
\newcommand{\bdm}{\begin{displaymath}}
        \newcommand{\edm}{\end{displaymath}}
\newcommand{\bds} {\begin{description}}
        \newcommand{\eds} {\end{description}}
\newcommand{\ben}{\begin{enumerate}}
        \newcommand{\een}{\end{enumerate}}
\newcommand{\beq}{\begin{equation}}
        \newcommand{\eeq}{\end{equation}}
\newcommand{\bfg} {\begin{figure}[h]}
        \newcommand{\efg} {\end{figure}}
\newcommand{\bi} {\begin {itemize}}
        \newcommand{\ei} {\end {itemize}}
\newcommand{\bqn}{\begin{eqnarray}}
        \newcommand{\eqn}{\end{eqnarray}}
\newcommand{\bqs}{\begin{eqnarray*}}
        \newcommand{\eqs}{\end{eqnarray*}}
\newcommand{\bsl} {\begin{slide}[8.8in,6.7in]}
        \newcommand{\esl} {\end{slide}}
\newcommand{\bsq}{\begin{subequations}}
        \newcommand{\esq}{\end{subequations}}       
\newcommand{\bss} {\begin{slide*}[9.3in,6.7in]}
        \newcommand{\ess} {\end{slide*}}
\newcommand{\btb} {\begin {table}}
        \newcommand{\etb} {\end {table}}
\newcommand{\bvb}{\begin{verbatim}}
        \newcommand{\evb}{\end{verbatim}}

\newcommand{\m}{\mbox}

\newcommand {\pd}[2] {{\frac {\partial {#1}} {\partial {#2}}}}

\newcommand{\cas}[1]{{{\left \{ \ba #1 \ea \right. }}}

\newcommand{\reff}[1] {{{Figure \ref {#1}}}}
\newcommand{\refe}[1] {{(\ref {#1})}}

\def\pmb#1{\setbox0=\hbox{$#1$}%
   \kern-.025em\copy0\kern-\wd0
   \kern.05em\copy0\kern-\wd0
   \kern-.025em\raise.0433em\box0 }

\def\eop{{\hfill $\blacksquare$}}

\newtheorem{theorem}{Theorem}[section]
\newtheorem{definition}[theorem]{Definition}

\newtheorem{lemma}[theorem]{Lemma}

\def\B {{\mathbb{B}}}

\usepackage[draft]{changes}
\definechangesauthor[name={Wenlong Jin},color=red]{WJ}

\usepackage{multirow}
\usepackage{subcaption}
\usepackage{float}
\usepackage[colorinlistoftodos]{todonotes}

\begin{document}
	\title{On the global stability of departure time user equilibrium: A Lyapunov approach} 
\author{Wen-Long Jin \footnote{Department of Civil and Environmental Engineering, California Institute for Telecommunications and Information Technology, Institute of Transportation Studies, 4000 Anteater Instruction and Research Bldg, University of California, Irvine, CA 92697-3600. Tel: 949-824-1672. Fax: 949-824-8385. Email: wjin@uci.edu. Corresponding author}}
\maketitle
\begin{abstract}
In \citep{jin2018_stable}, a new day-to-day dynamical system was proposed for drivers' departure time choice at a single bottleneck. Based on three behavioral principles, including backward choice, cost balancing, and scheduling cost reducing or scheduling payoff gaining principles, the nonlocal departure and arrival times choice problems were converted to the local scheduling payoff choice problem, whose day-to-day dynamics are described by the Lighthill-Whitham-Richards (LWR) model on an imaginary road of increasing scheduling payoff. Thus the departure time user equilibrium (DTUE), the arrival time user equilibrium (ATUE), and the scheduling payoff user equilibrium (SPUE) are uniquely determined by the stationary state of the LWR model, which was shown to be locally, asymptotically stable with analysis of the discrete approximation of the LWR model and through a numerical example.

In this study attempt to analytically prove the global stability of the SPUE, ATUE, and DTUE.
We first generalize the conceptual models for arrival time and scheduling payoff choices developed in \citep{jin2018_stable} for a single bottleneck with a generalized scheduling cost function, which includes the cost of the free-flow travel time. Then we present the LWR model for the day-to-day dynamics for the scheduling payoff choice as well as the SPUE. 
We further formulate a new optimization problem for the SPUE and demonstrate its equivalent to the optimization problem for the ATUE in \citep{iryo2007equivalent}. Finally we show that the objective functions in the two optimization formulations are equal and can be used as the potential function for the LWR model and prove that the stationary state of the LWR model, and therefore, the SPUE, DTUE, and ATUE, are globally, asymptotically stable, by using Lyapunov's second method. Such a globally stable behavioral model can provide more efficient departure time and route choice guidance for human drivers and connected and autonomous vehicles in more complicated networks.	

\end{abstract}
 
\section{Introduction}
If an origin-destination pair is connected by a single bottleneck, morning commuters would diversify their departure times so as to minimize their total costs, including both travel costs, scheduling costs, and/or tolls \citep{vickrey1969congestion,hendrickson1981schedule,small2015bottleneck}. Theoretically, the traffic system reaches the so-called departure time user equilibrium (DTUE), when ``the journey costs at all  departure times actually used are equal, and (equal to or) less than those which would be experienced by a single vehicle at any unused time'' \citep{wardrop1952ue,Arnott1990bottleneck}.  From the point view of arrival times, the DTUE is equivalent to the arrival time user equilibrium (ATUE), where ``the journey costs at all arrival times actually used are equal, and (equal to or) less than those which would be experienced by a single vehicle at any unused time''. In the ATUE, all commuters have the same total cost, and ``no individual user can improve his total cost by unilaterally changing arrival times" \citep{hendrickson1981schedule,mahmassani1984due}.\footnote{In this study the departure and arrival times of a vehicle represent the time for it to depart from the origin and arrive at the destination.}

Traditionally, the DTUE and ATUE have been analyzed and solved through variational inequalities \citep[e.g.,][]{Friesz1993due,Szeto2004departure}, optimization (linear programming) \citep{iryo2007equivalent}, and linear complementarity formulations \citep{akamatsu2015corridor}. 
In contrast to such purely phenomenological approaches, day-to-day dynamical system models attempt to explain why and how a DTUE/ATUE can be reached through choice and learning behaviors. Many such models have been successfully applied to study route choice behaviors and shown to be stable and converge to the corresponding static user equilibrium \citep{smith1984stability,Friesz1994dynamic,Nagurney1996ds,Nagurney1997equilibria,jin2007dstap,yang2009day,xiao2016physics,guo2016discrete}.  All these models are based on the fundamental behavioral principle that drivers tend to switch to less costly routes, but different in their implementation details. However, simple extensions of this behavioral principle for departure time choice have been shown to be unsuccessful \citep{iryo2008analysis,bressan2012variational,guo2016departure,guo2017day}; that is, all the proposed dynamical systems in these references are unstable. The lack of a stable day-to-day dynamical system model has led to questions over the existence of a stable DTUE/ATUE in the real world.

In \citep{jin2018_stable}, it was argued that the existence of a stable day-to-day dynamical system for departure time choice cannot be ruled out logically, and observations of relatively stationary day-to-day traffic patterns as well as relatively fixed departure time choices of commuters suggest stable day-to-day dynamics for departure time choice. However, the behavioral principles regarding departure time choice should be more sophisticated than directly switching to a less costly departure or arrival time.
  Further three behavioral principles were identified: (i) the backward choice principle: drivers choose their arrival times before departure times; (ii) the cost balancing principle for departure time choice: drivers choose their departure times to balance the total costs; and (iii) the scheduling cost reducing or scheduling payoff gaining principle for arrival time choice: drivers switch their arrival times to those when the bottleneck is under-utilized with larger scheduling payoffs. Therefore, even though the DTUE and ATUE are equivalent and seemingly symmetric, drivers' departure time and arrival time choices are not, and the former is dictated by the latter. Further with a conceptual V-shaped tube model, it was shown that the arrival time choice is nonlocal, with the set of target arrival times disconnected. However, by converting the V-shaped tube into a single tube, one can show that the nonlocal arrival times choice problem is equivalent to the local scheduling payoff choice problem. Thus the day-to-day departure time choice leads to traffic flow evolution on an imaginary road, which can be naturally described by the Lighthill-Whitham-Richards (LWR). It was further proved that the DTUE and ATUE are equivalent to the scheduling payoff user equilibrium (SPUE) defined by the LWR model and the corresponding splitting and cost balancing procedures to determine the equilibrated arrival and departure flow-rates. Such a new day-to-day dynamical system was shown to be asymptotically stable for the discrete approximation of the LWR model analytically and then through a numerical example. However, the analytical proof is for local stability subject to small perturbations, and the numerical proof is just for one example.

In this study we attempt to analytically prove the global stability of the SPUE, ATUE, and DTUE.
We first generalize the conceptual models for arrival time and scheduling payoff choices developed in \citep{jin2018_stable} for a single bottleneck with a generalized scheduling cost function, which includes the cost of the free-flow travel time. Then we present the LWR model for the day-to-day dynamics for the scheduling payoff choice as well as the SPUE. 
We further formulate a new optimization problem for the SPUE and demonstrate its equivalent to the optimization problem for the ATUE in \citep{iryo2007equivalent}. Finally we show that the objective functions in the two optimization formulations can be used as the potential function for the LWR model and prove that the SPUE and, therefore, the DTUE and ATUE, are asymptotically stable, by using Lyapunov's second method.

The rest of the paper is organized as follows. In Section 2, we present the definitions of the single bottleneck problem and conceptual models for day-to-day arrival time choice. In Section 3, we present the dynamical system model of scheduling payoff choice and discuss the SPUE.  In Section 4, we formulate an optimization problem for the SPUE and demonstrate its equivalent to the optimization problem for the ATUE in \citep{iryo2007equivalent}. In Section 5, we show that the objective functions in the two optimization formulations can be used as the potential function for the LWR model and prove that the SPUE and, therefore, the DTUE and ATUE, are asymptotically stable, by using Lyapunov's second method. In Section 6, we conclude the study with discussions and future directions.

\section{Definitions and conceptual models}

\subsection{Definitions}
We consider a road with a single bottleneck between an origin and destination pair, where the free-flow travel time, $\Upsilon_0$, is constant from day to day. We assume that  the bottleneck capacity (maximum service rate), $C$, and the travel demand, i.e., the total number of vehicles, $N$, are also constant. 
As illustrated in \reff{bottleneck_traffic_costs}(a), on day $r$, the departure and arrival cumulative flows at time $t$ are respectively denoted by $F'(r,t)$ (thick, red curve) and $G(r,t)$ (thin, blue curve). We shift $F'(r,t)$ to the right by $\Upsilon_0$ and obtain $F(r,t)=F'(r,t-\Upsilon_0)$. The derivatives of $F'(r,t)$, $F(r,t)$, and $G(r,t)$ are respectively denoted by $f'(r,t)$, $f(r,t)$, and $g(r,t)$; thus $f'(r,t)=f(r,t+\Upsilon_0)$ and $g(r,t)$ are respectively the departure and arrival flow-rates. Since the departure time choice is dictated by the arrival time choice, we are primarily concerned with the travel characteristics for  vehicles arriving at $t$:  $\delta(r,t)$ is the queue size, and $\Upsilon(r,t)$ the queueing time. Thus the travel time for vehicles arriving at the destination at $t$ is $\Upsilon_0+\Upsilon(r,t)$. For a vehicle with a departure time of $t$, its arrival time is $t+\Upsilon_0+\frac{\delta(r,t+\Upsilon_0)}C$. If $g(r,t)<C$, the bottleneck is under-utilized at $t$; otherwise $g(r,t)=C$, and the bottleneck is fully utilized. The queue length is zero when the bottleneck is under-utilized.

The total cost for vehicles arriving at the destination at $t$ is denoted by $\phi(r,t)$, which comprises of the travel cost, caused by the travel time, and the scheduling cost, caused by the schedule delay:
\bqn
\phi(r,t)&=&\alpha (\Upsilon_0+\Upsilon(r,t))+\beta \max\{ t_*-t,0\} +\gamma \max\{t-t_*, 0\}, \label{def:totalcost}
\eqn
where $t_*$ is the ideal arrival time.
In \citep{arnott1990departure}, $\alpha=$\$6.4/hr, $\beta=$\$3.90/hr, and $\gamma=$15.21/hr. A necessary condition for the existence of DTUE and ATUE is that $\beta<\alpha$ \citep{small2015bottleneck}.
Further we denote the queueing cost by 
\bqn
\phi_1(r,t)&=&\alpha \Upsilon(r,t), \label{def:cost1}
\eqn
and the non-queueing cost, including both the free-flow travel cost and the scheduling cost, by
\bqn
\phi_2(t)&=&\alpha \Upsilon_0+\beta \max\{ t_*-t,0\} +\gamma \max\{t-t_*, 0\}. \label{def:cost2}
\eqn
We refer to $\phi_2(t)$ as generalized scheduling cost.
Here we assume that the coefficients, $\alpha$, $\beta$, and $\gamma$, are constant from day to day. Thus the queueing cost may vary from day to day but the other cost does not. \reff{bottleneck_traffic_costs}(b) illustrates the cost functions for given departure and arrival cumulative flows. Further we have the following theorem from \citep{jin2018_stable}.

\begin{theorem} \label{theorem:vacancy} 
	If the bottleneck is under-utilized at $t_1$ and the generalized scheduling cost at $t_2$ is not smaller than that at $t_1$, then the total cost at $t_2$ is not smaller than that at $t_1$. That is, if $g(r,t_1)<C$, and $\phi_2(t_2)\geq\phi_2(t_1)$, then $\phi(r,t_2) \geq \phi(r,t_1)$.
\end{theorem}

		\bfg\bc
		\begin{tabular}{cc}
			\includegraphics[width=2.5in]{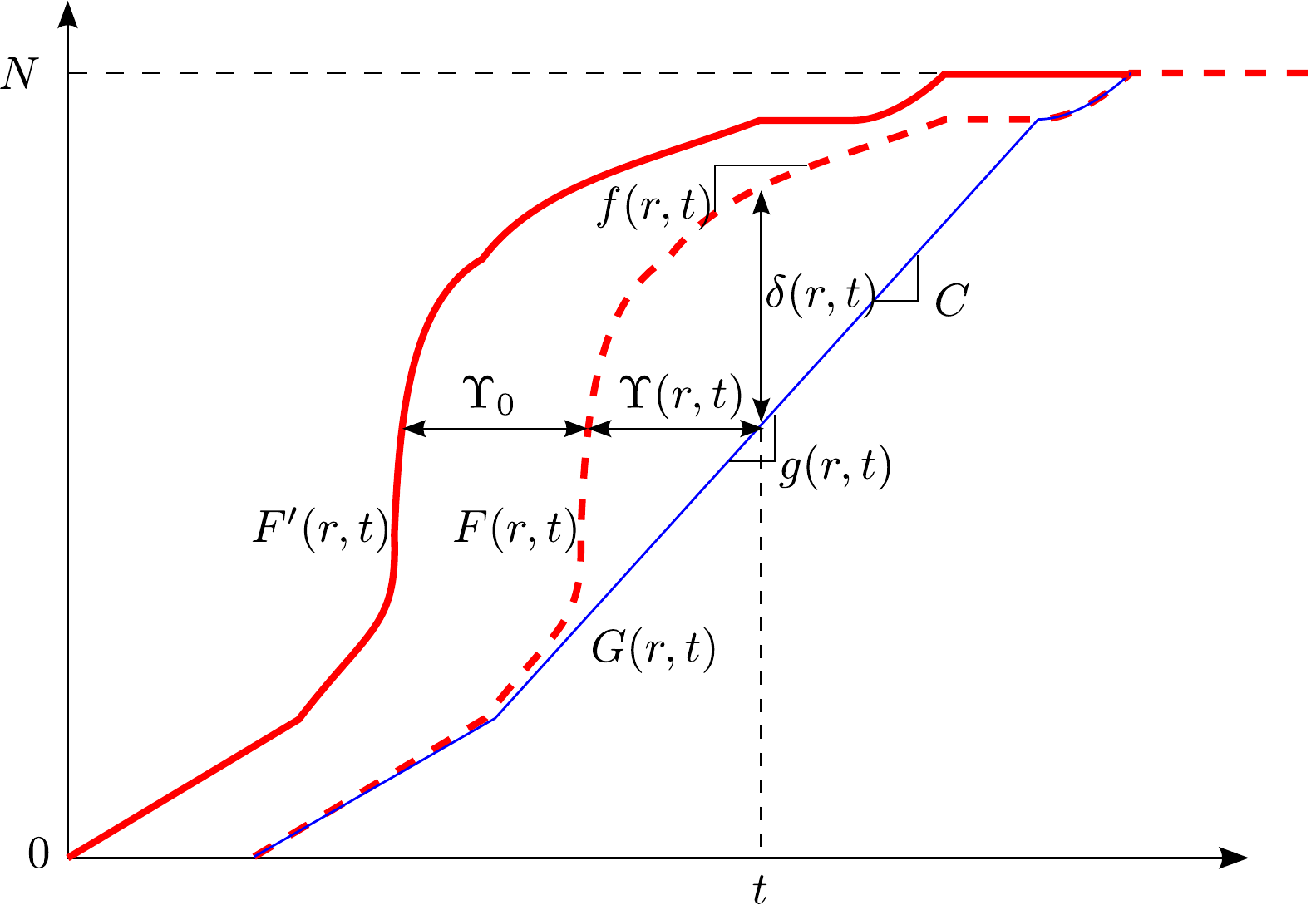} & \includegraphics[width=3.5in]{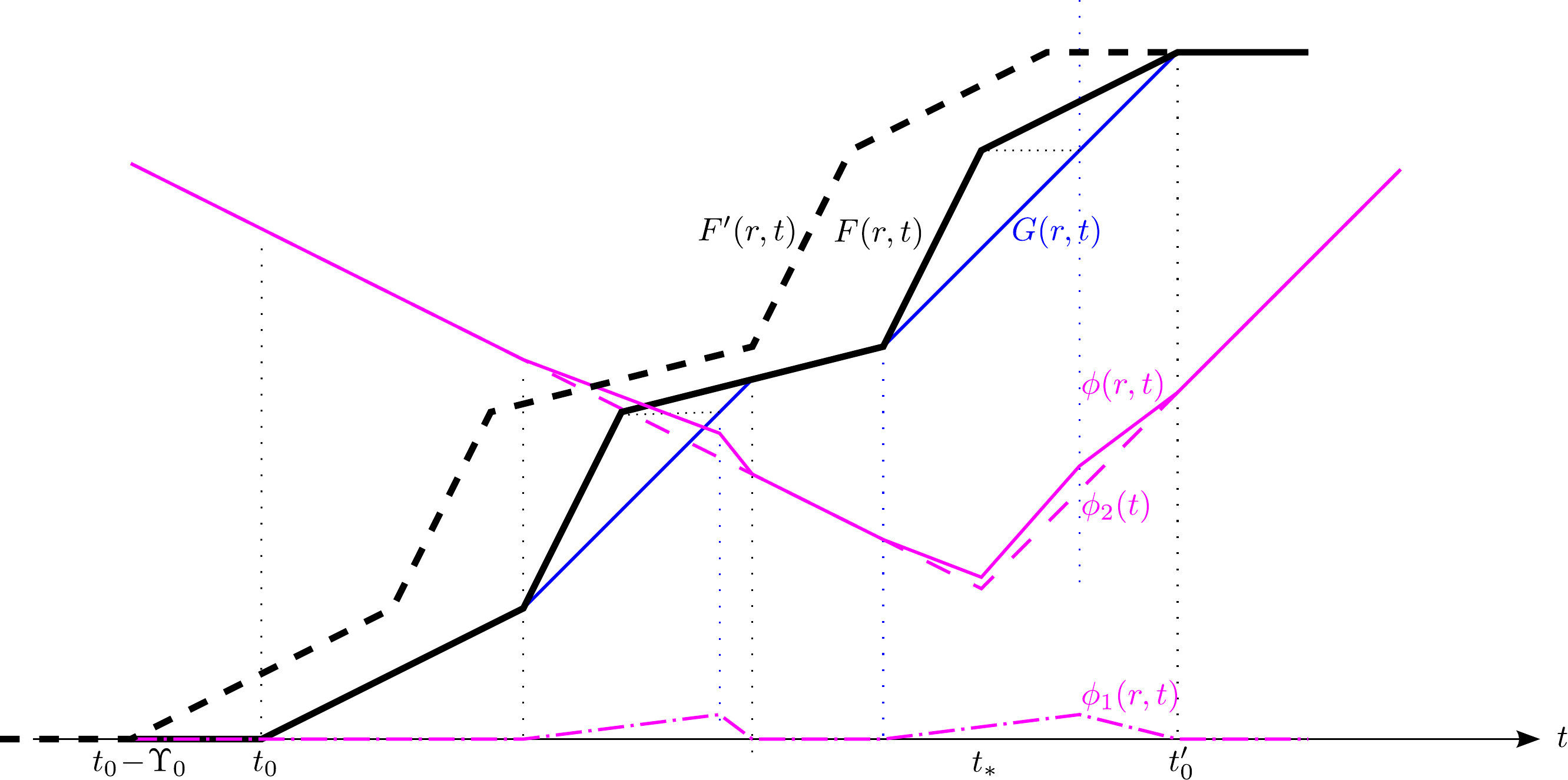} \\
			(a) & (b) \\
		\end{tabular}
		\caption{(a) Traffic variables for a single bottleneck  on day $r$; (b) Cost functions} \label{bottleneck_traffic_costs}
		\ec\efg

\begin{definition} \label{def:DTUE:ATUE}
	The system reaches the {\bf arrival time user equilibrium (ATUE)} if all used arrival times have the same total cost, which is not greater than those of unused arrival times. That is, $\phi(r,t)=\phi^*$ if $g(r,t)>0$; and $\phi(r,t)\geq \phi^*$ if $g(r,t)=0$. Here $\phi^*$ is the minimum total cost.
	
	The system reaches the {\bf departure time user equilibrium (DTUE)} if all used departure times have the same total cost, which is not greater than those of unused departure times. That is, $\phi\left(r,t+\Upsilon_0+\frac{\delta(r,t+\Upsilon_0)}C\right)=\phi^*$ if $f'(r,t)>0$; and $\phi\left(r,t+\Upsilon_0+\frac{\delta(r,t+\Upsilon_0)}C\right)\geq \phi^*$ if $f'(r,t)=0$. 
\end{definition}

The equivalence between ATUE and DTUE is apparent, since each commuter has a unique set of arrival and departure times. A mathematical proof of the equivalence was presented in \citep{jin2018_stable}. Since a vehicle's departure time choice is dictated by its arrival time choice, we focus on finding the ATUE.

\subsection{Conceptual models for day-to-day arrival time choice}

\bfg\bc
\includegraphics[width=5.5in]{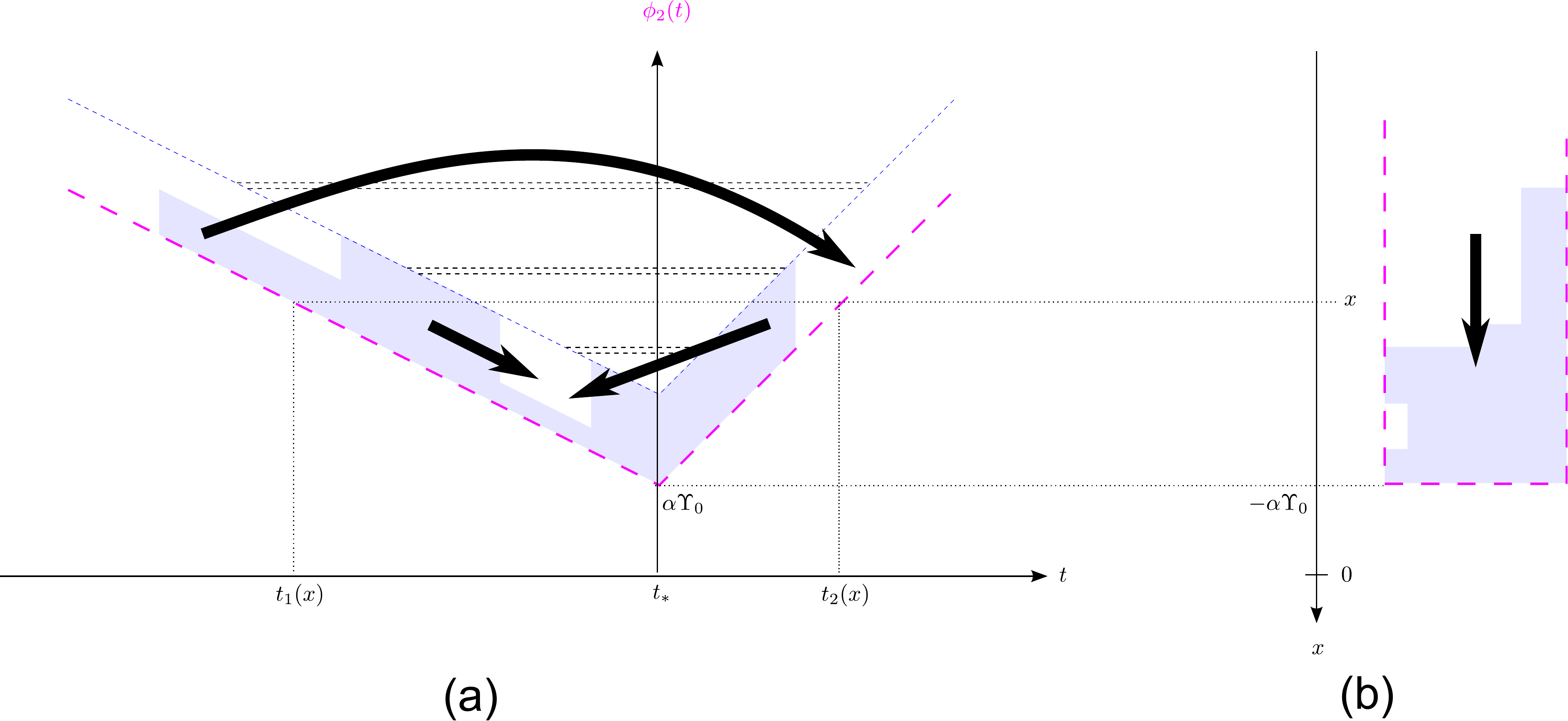}\caption{Conceptual models of arrival time and scheduling payoff choices} \label{single_tube_general}
\ec\efg

In \reff{single_tube_general}(a),  at $t$ the height of the dashed curve represents the generalized scheduling cost, $\phi_2(t)$,  the vertical gap between the dashed and dotted curves represents the bottleneck capacity, $C$, and the height of the shaded region represents the arrival flow-rate, $g(r,t)$. In this case, the region bounded by the dashed and dotted curves forms a V-shaped tube, which is under-utilized at $t$ if  there is vacancy in the tube. According to Theorem \ref{theorem:vacancy}, vehicles can reduce their total costs by switching to under-utilized arrival times with smaller generalized scheduling costs; this is the scheduling cost reducing principle for arrival time choice. In the figure, the three arrows represent three of such switches of arrival times. 
Conceptually such choices are consistent with the movements of fluids (e.g. water) caused by gravity in a V-shaped tube: they will attempt to fill a vacancy at a lower point. In this sense, the third principle of arrival time choice is equivalent to finding a lower vacancy  in the V-shaped tube.  

However,  since drivers can switch to  arrival times earlier or later than the ideal arrival time, the fluid movements can be nonlocal \citep{bressan2012variational}, which is enabled if the two sides of the V-shaped tube are connected by capillaries with negligible widths. Thus, the V-shaped tube is equivalent to a single tube as shown in \reff{single_tube_general}(b), where the fluid movements are local. The single tube can also be imagined as a road, and vehicles drive southbound on the imaginary road.

We introduce a new variable $x$ for the scheduling payoff, equal to the negative generalized scheduling cost, as shown in \reff{single_tube_general}(b): 
\bqn
x&=&-\phi_2(t) \leq -\alpha \Upsilon_0. \label{def:x}
\eqn
Therefore, a value of $x$ corresponds to two arrival times, $t_1(x)$ and $t_2(x)$, as shown in \reff{single_tube_general}(a), where
\bsq\label{def:t1t2}
\bqn
t_1(x)&=&t_*+\frac{x+\alpha \Upsilon_0}{\beta}, \\
t_2(x)&=&t_*-\frac{x+\alpha \Upsilon_0}{\gamma},
\eqn
\esq
where $x+\alpha \Upsilon_0\leq 0$ and $t_1(x)\leq t_*\leq t_2(x)$. Thus vehicles drives on the imaginary road in the positive direction of $x$; i.e., the increasing direction of the scheduling payoff.

\section{A dynamical system model and the scheduling payoff user equilibrium}
\subsection{The LWR model for day-to-day scheduling payoff choice}
On the imaginary road, the density at the scheduling payoff of $x$ is denoted by $k(r,x)$, whose unit is veh/\$ and defined by the integral form: 
\bqn
\int_{y=x}^{-\alpha \Upsilon_0} k(r,y) dy&=& \int_{t_1(x)}^{t_2(x)} g(r,t) dt, \label{def:rho_int}
\eqn
or the differential form: 
\bqn
k(r,x)&=& \frac 1 \beta g(r,t_1(x))  + \frac 1 \gamma g(r,t_2(x)). \label{def:rho}
\eqn
Therefore the imaginary density at $x$ equals the sum of the arrival flow-rates at the two time, $t_1(x)$ and $t_2(x)$, adjusted by $\beta$ and $\gamma$ respectively. It is represented by the horizontal width at $x$ of the shaded region in \reff{single_tube_general}(b).
Correspondingly, the width of the single tube at $x$ shown in \reff{single_tube_general}(b) equals the sum of the capacities at the two times, adjusted by $\beta$ and $\gamma$ respectively. We refer to the width as the ``jam density", which is denoted by $\kappa$; thus
\bqn
\kappa=(\frac 1\beta+\frac 1\gamma)C. \label{def:kappa}
\eqn

Furthermore, from
\bqn
\int_{x=-\infty}^{-\alpha \Upsilon_0} k(r,x) dx&=&\int_{t=-\infty}^{\infty} g(r,t) dt=N, \label{total_conservation}
\eqn
we can see that the imaginary density is conserved, and 
\bqn
\pd{}r k(r,x)+ \pd{}x q(r,x)&=&0. \label{lwr_conservation}
\eqn
The number of vehicles switching from payoff $x^-$ to $x^+$ on day $r$ is denoted by $q(r,x)$, which is the day-to-day flow-rate on the imaginary road with a unit of veh/day. 
We assume that there exists a fundamental diagram:
\bqn
q&=&Q(k). \label{lwr_fd}
\eqn
 \footnote{
An example is the triangular fundamental diagram:	
	\bqn
	Q(k)&=&\min\{u k, w(\kappa-k) \}, \label{lwr_tri_fd}
	\eqn
	where $u$ is the free-flow speed and $w$ the shock wave speed in congested traffic. Here both $u$ and $w$ are positive and can choose arbitrary values with a unit of \$/day.
But the dynamical system model works for any reasonable fundamental diagrams if $Q(k)\geq 0$, and $Q(k)=0$ if and only if $k=0$ or $k=\kappa$.}

Therefore, the traffic flow dynamics on the imaginary road can be naturally described by the Lighthill-Whitham-Richards (LWR) model  \citep{lighthill1955lwr,richards1956lwr}:
\bqn
\pd{}r k(r,x)+ \pd{}x Q(k(r,x))&=&0, \label{lwr_sc}
\eqn
which is derived from \refe{lwr_conservation} and \refe{lwr_fd}.

In addition, the LWR model satisfies the following conditions:
\ben
\item $k(r,x) \in[0,\kappa]$; 
\item The initial condition:
\bqn
k(0,x)&=&k_0(x), \quad x\in (-\infty, -\alpha \Upsilon_0],\label{lwr_ini}
\eqn
which can be calculated from the initial departure cumulative flows $F'(0,t)$ and the point queue model;
\item The boundary conditions: 
\bqn
q(r,-\alpha \Upsilon_0)&=&q(r,-\infty)=0; \label{lwr_bc}
\eqn
i.e., no vehicles can enter or leave the imaginary road.
\een

The LWR model, \refe{lwr_sc}, subject to the initial and boundary conditions, \refe{lwr_ini} and \refe{lwr_bc}, is an infinite-dimensional dynamical system model, in which the imaginary densities at different payoffs evolve with respect to the day variable, $r$. Thus we convert the departure time choice problem into a day-to-day traffic flow problem on an imaginary road. Note that the arrival and departure flow-rates can be calculated from the imaginary density by following the splitting and cost balancing procedures, which were presented in \citep{jin2018_stable} but omitted in this study.

\subsection{Scheduling payoff user equilibrium}

\begin{definition}
The system reaches the {\bf scheduling payoff user equilibrium (SPUE)} if for the LWR model, \refe{lwr_sc}, reaches a stationary state \citep{jin2012statics,jin2015existence}:
\bqn
\pd{k(r,x)}r &=&0, \label{def:SPUE}
\eqn
and the total costs are balanced by the cost balancing principle for departure time choice. Thus in the SPUE,
\bqn
k(r,x)&=&k^*(x).
\eqn
\end{definition}

\begin{theorem} The SPUE is equivalent to the ATUE and DTUE.
	\end{theorem}

This theorem was proved in  \citep{jin2018_stable}. Thus the LWR model,  \refe{lwr_sc}, is also the dynamical system model of day-to-day departure time choice, and ATUE/DTUE/SPUE is the stationary state of the LWR model.

\begin{lemma} \label{lemma_spue}In the SPUE, $q(r,x)=0$; i.e., there is no scheduling payoff choice dynamics.
	\end{lemma}
{\em Proof}. From \refe{lwr_conservation}, we have
\bqs
\pd{q(r,x)}x&=&0.
\eqs 
From \refe{lwr_fd} and \refe{def:SPUE}, we have 
\bqs
\pd{q(r,x)}r&=&0.
\eqs
Thus $q(r,x)$ is independent of both $r$ and $x$ and, therefore, constant.

Further from the boundary condition, \refe{lwr_bc}, we conclude that $q(r,x)=0$ in the SPUE.
\eop

\begin{theorem} In the SPUE, the density is uniquely given by
	\bqn
	k^*(x)&=&\cas{{ll} \kappa, & -\alpha \Upsilon_0- L^* \leq x\leq -\alpha \Upsilon_0; \\ 0, &\m{otherwise},} \label{stationarydensity}
	\eqn
	where
	\bqn
	L^*&=&\frac{N} {\kappa}. \label{def:Lstar}
	\eqn
	That is, all vehicles have their generalized scheduling costs not greater than $\alpha \Upsilon_0+L^*$, which is the total cost for all vehicles, and the corresponding smallest and largest arrival times are $t_1(- \alpha \Upsilon_0-L^*)$ and $t_2(-\alpha \Upsilon_0-L^*)$.
This can be illustrated in \reff{single_tube_general}(b), where the fluid settles down to the bottom part of the single tube when there is no scheduling payoff choice dynamics.
	\end{theorem}
{\em Proof}. From Lemma \ref{lemma_spue}, $q(r,x)=0$, and $k^*(x)=0$ or $\kappa$ at any $x$. Further from the kinematic wave theory \citep{jin2009sd} we have
\bqn
q(r,x)&=&\min\{d(r,x^-), s(r,x^+)\}=\min\{Q( \min \{k(r,x^-), \kappa_c \}), Q( \max \{k(r,x^+), \kappa_c \}) \}, \label{lwr_kwt}
\eqn
where $\kappa_c$ is the critical density, and $d(r,x^-)$ and $s(r,x^+)$ are respectively the upstream demand and downstream supply. Thus we can have the following three possibilities at any $x$: (i) $k(r,x^-)=k(r,x^+)=0$; (ii)  $k(r,x^-)=k(r,x^+)=\kappa$; and (iii) $k(r,x^-)=0$, and $k(r,x^+)=\kappa$.  In particular, $k(r,x^-)=\kappa$, and $k(r,x^+)=0$ are not allowed. Therefore, inside the single tube or imaginary road, the stationary density has to be given by \refe{stationarydensity}, which constitutes a zero-speed shock wave at $-\alpha \Upsilon_0-L^*$.

Furthermore from \refe{total_conservation}, we can find $L^*$ in \refe{def:Lstar}, which is the length of the queue. \eop

\section{Two optimization formulations}

\subsection{A new optimization formulation of the SPUE}
We define the following functional of $k(r,x)$:
\bqn
\Phi(k(r,x))&=&- \int_{-\infty}^{-\alpha \Upsilon_0} x k(r,x) dx. \label{def:potential}
\eqn

\begin{lemma} \label{minimization_lemma}
	$\Phi(k(r,x))$ reaches its minimum if and only if $k(r,x)=k^*(x)$ given by \refe{stationarydensity}, where
	\bqn
	\Phi(k^*(x))&=&- \kappa \int_{-\alpha \Upsilon_0-L^*}^{-\alpha \Upsilon_0} x  dx.
	\eqn
\end{lemma}
{\em Proof}. If $k(r,x)\neq k^*(x)$, then there exists $x<-\alpha\Upsilon_0-L^*$ such that $k(r,x)>0$, and $\int_{-\infty}^{-\alpha\Upsilon_0-L^*} k(r,x)dx=B>0$.\footnote{Note that when $k(r,x)>0$ for some $x<-\alpha\Upsilon_0-L^*$, and $\int_{-\infty}^{-\alpha\Upsilon_0-L^*} k(r,x) dx=0$, then $k(r,x)$ and $k^*(x)$ are only different on a set of measure 0. Such a difference is not physically meaningful, and we exclude this case in our study.}  Therefore,
\bqs
\frac 1 B [\Phi(k(r,x))-\Phi(k^*(x))]&=&-\int_{-\infty}^{-\alpha \Upsilon_0-L^*} x \frac {k(r,x)}B dx +  \int_{-\alpha \Upsilon_0-L^*}^{-\alpha \Upsilon_0} x  \frac{\kappa-k(r,x)}B dx.
\eqs

From \refe{total_conservation} we have
\bqs
\int_{-\infty}^{-\alpha \Upsilon_0}  k(r,x) dx&=&N=\int_{-\alpha \Upsilon_0-L^*}^{-\alpha \Upsilon_0}  \kappa dx,
\eqs
which leads to
\bqs
\int_{-\infty}^{-\alpha \Upsilon_0-L^*}  \frac{k(r,x)}B dx&=&\int_{-\alpha \Upsilon_0-L^*}^{-\alpha \Upsilon_0}  \frac{\kappa-k(r,x)}B dx=1.
\eqs
Since $0\leq k(r,x)\leq \kappa$, both $\frac{k(r,x)}B$ and $\frac{\kappa-k(r,x)}B$ are non-negative. Therefore, $\int_{-\infty}^{-\alpha \Upsilon_0-L^*} x \frac {k(r,x)}B dx$ is an average of $x$ between $-\infty$ and $-\alpha \Upsilon_0-L^*$. 
Note that there exists $x<-\alpha\Upsilon_0-L^*$ such that $k(r,x)>0$. Thus we have
\bqs
\int_{-\infty}^{-\alpha \Upsilon_0-L^*} x \frac {k(r,x)}B dx<-\alpha \Upsilon_0-L^*.
\eqs
Similarly, $\int_{-\alpha \Upsilon_0-L^*}^{-\alpha \Upsilon_0} x  \frac{\kappa-k(r,x)}B dx$ is an average of $x$ between $-\alpha \Upsilon_0-L^*$ and $-\alpha \Upsilon_0$, and
\bqs
\int_{-\alpha \Upsilon_0-L^*}^{-\alpha \Upsilon_0} x  \frac{\kappa-k(r,x)}B dx\geq -\alpha \Upsilon_0-L^*.
\eqs

Therefore $\frac 1 B [\Phi(k(r,x))-\Phi(k^*(x))]>0$, which leads to $\Phi(k(r,x))>\Phi(k^*(x))$ when $k(r,x)\neq k^*(x)$, and $\Phi(k(r,x))$ reaches its minimum if and only if $k(r,x)=k^*(x)$. 
\eop

From Lemma \ref{minimization_lemma}, we have the following theorem.

\begin{theorem} \label{minimization_theorem}
	$k^*(x)$ in \refe{stationarydensity} is the unique solution of the following optimization problem:
\bqn
\min_{k(r,x)} \Phi(k(r,x))&=&- \int_{-\infty}^{-\alpha \Upsilon_0} x k(r,x) dx
\eqn
s.t.
\bsq\label{k_constraints}
\bqn
\int_{-\infty}^{-\alpha \Upsilon_0} k(r,x) dx&=&N,\\
k(r,x)&\geq&0,\\
k(r,x)&\leq&\kappa.
\eqn
\esq
This is a new optimization formulation of the SPUE/ATUE/DTUE and the corresponding choice behaviors.
\end{theorem}

\subsection{An existing optimization formulation of the ATUE}
In the ATUE, the equilibrium arrival flow-rates corresponding to $k^*(x)$ in \refe{stationarydensity} are given by
\bqn
g^*(t)&=&\cas{{ll} C, & t\in\B'; \\
	0, &t\in \B,} \label{equilibrium_arrival_flowrate}
\eqn
where the interval $\B=(-\infty, t_1(-\alpha \Upsilon_0 -L^*))\cup(t_2(-\alpha \Upsilon_0 -L^*),\infty)$ and its complement  $\B'=[t_1(-\alpha \Upsilon_0 -L^*),t_2(-\alpha \Upsilon_0 -L^*)]$.

In \citep{iryo2007equivalent}, a discrete optimization formulation of the ATUE was defined for a single bottleneck with multiple classes. The continuous version for a single class can be written as follows.

\begin{theorem}  The equilibrium arrival flow-rates, $g^*(t)$, is the unique solution of the following optimization problem:
	\bqn
	\min_{g(r,t)}\Phi'(g(r,t))&=&\int_{-\infty}^{\infty} g(r,t) \phi_2(t) dt, \label{def:ATUEobjective}
	\eqn
	s.t.
	\bqn
	\int_{-\infty}^\infty g(r,t)dt&=&N,\label{total_g_constraint}\\
	g(r,t)&\geq&0,\\
	g(r,t)&\leq&C.
	\eqn
	In particular
	\bqn
	\Phi'(g^*(t))&=&\int_{t\in\B'} C \phi_2(t) dt.
	\eqn
	\end{theorem}
{\em Proof}. First, $g^*(t) \geq 0$, $g^*(t)\leq C$, and from \refe{def:t1t2}, \refe{def:kappa}, and \refe{def:Lstar},  we have
\bqs
\int_{-\infty}^\infty g^*(t)dt&=&C \cdot (t_2(-\alpha \Upsilon_0 -L^*) -t_1(-\alpha \Upsilon_0 -L^*))=C L^* (\frac 1\beta+\frac 1\gamma)=N.
\eqs
Thus $g^*(t)$ satisfies the constraints.

If a feasible $g(r,t)\neq g^*(t)$, then there exists $t\in\B$ such that $g(r,t)>0$, and $\int_{t\in\B} g(r,t)dt=B>0$. \footnote{We exclude the case when $g(r,t)>0$ for some $t\in\B$ , but $\int_{t\in\B} g(r,t)dt=0$, since such a $g(r,t)$ is only different from $g^*(t)$ on a set of measure 0 and not physically meaningful.}
Therefore,
\bqs
\frac 1B[\Phi'(g(r,t))-\Phi'(g^*(t))]&=&\int_{t\in\B} \phi_2(t) \frac{g(r,t)}B dt- \int_{t\in \B'} \phi_2(t) \frac{C-g(r,t)}B dt.
\eqs

From \refe{total_g_constraint} we have 
\bqs
\int_{t\in \B} g(r,t)dt+\int_{t\in\B'} g(r,t)dt&=&N=\int_{t\in \B'} g^*(t) dt.
\eqs
Thus
\bqs
\int_{t\in\B} \frac{g(r,t)}B dt&=&\int_{t\in \B'}  \frac{C-g(r,t)}B dt=1.
\eqs
Since $0\leq g(r,t)\leq C$, both $\frac{g(r,t)}B$ and $\frac{C-g(r,t)}B$ are non-negative. Therefore, $\int_{t\in\B} \phi_2(t) \frac{g(r,t)}B dt$ is an average of $\phi_2(t)$ inside $\B$, and 
\bqs
\int_{t\in\B} \phi_2(t) \frac{g(r,t)}B dt&>&\phi_2(t_1(-\alpha\Upsilon_0-L^*))=\phi_2(t_2(-\alpha\Upsilon_0-L^*)).
\eqs
Similarly, $\int_{t\in \B'} \phi_2(t) \frac{C-g(r,t)}B dt$ is an average of $\phi_2(t)$ inside $\B'$, and 
\bqs
\int_{t\in \B'} \phi_2(t) \frac{C-g(r,t)}B dt&\leq&\phi_2(t_1(-\alpha\Upsilon_0-L^*))=\phi_2(t_2(-\alpha\Upsilon_0-L^*)).
\eqs

Therefore, $\frac 1B[\Phi'(g(r,t))-\Phi'(g^*(t))]>0$, which leads to $\Phi'(g(r,t))>\Phi(g^*(t))$ when $g(r,t)\neq g^*(t)$. Hence $g^*(t)$ solves the minimization problem. \eop

\subsection{Equivalence of the two optimization formulations}

\begin{theorem}
	Given \refe{def:x} and \refe{def:rho}, the objective functions of the SPUE and ATUE optimization formulations are equal; i.e.,
	\bqn
	\Phi(k(r,x))&=&\Phi'(g(r,t)).
	\eqn
	The constraints are also equivalent.
	Therefore, the two optimization formulations are equivalent.
	\end{theorem}
{\em Proof}.
From \refe{def:x} and the relationship between $t$ and $x$ in \refe{def:t1t2}, we have that 
\bqs
\int_{-\infty}^{t_*} g(r,t) \phi_2(t) dt &=& \int_{-\infty}^{t_*} g(r,t) (\alpha \Upsilon_0+\beta(t_*-t)) dt =- \int_{-\infty}^{-\alpha\Upsilon_0} x g(r,t_1(x)) \frac 1\beta dx,
\eqs
and
\bqs
\int_{t_*}^{\infty} g(r,t) \phi_2(t)  dt&=&\int_{t_*}^{\infty} g(r,t) (\alpha \Upsilon_0+\gamma(t-t_*)) dt =\int_{-\alpha\Upsilon_0}^{-\infty} x g(r,t_2(x)) \frac 1\gamma dx.
\eqs

Thus we have
\bqs
\Phi'(g(r,t))&=&\int_{-\infty}^{t_*} g(r,t) \phi_2(t) dt +\int_{t_*}^{\infty} g(r,t) \phi_2(t)  d\\
&=&- \int_{-\infty}^{-\alpha\Upsilon_0} x (\frac 1\beta g(r,t_1(x)) +\frac 1\gamma g(r,t_2(x))) dx.
\eqs
Further from \refe{def:rho} we have
\bqs
\Phi'(g(r,t))&=&- \int_{-\infty}^{-\alpha\Upsilon_0} x k(r,x) dx=\Phi(k(r,x)).
\eqs
Thus the two objective functions are equal. It is easy to check that the constraints are also equivalent. Therefore the two optimization formulations are equivalent.
\eop

\section{Global stability via Lyapunov's second method}

\begin{lemma} \label{derivative_lemma}
	The day-derivative of $\Phi(k(r,x))$ is non-positive for any $k(r,x)$, and zero only at the SPUE, $k^*(x)$; i.e., for any feasible $k(r,x)$ given in \refe{k_constraints}, 
	\bqn
	\pd{\Phi(k(r,x))}r &=&0, \quad k(r,x) = k^*(x);\\
	\pd{\Phi(k(r,x))}r &<&0, \quad k(r,x) \neq k^*(x).	
	\eqn
	\end{lemma}
{\em Proof}. From \refe{def:potential} we have
\bqs
\pd{\Phi(k(r,x))}r&=&-\int_{-\infty}^{-\alpha\Upsilon_0} x \pd{k(r,x)} r dx.
\eqs
From the conservation equation, \refe{lwr_conservation}, we have $\pd{k(r,x)} r=-\pd{q(r,x)}x$, and
\bqs
\pd{\Phi(k(r,x))}r&=&\int_{-\infty}^{-\alpha\Upsilon_0} x \pd{q(r,x)} x dx=xq(r,x)|_{-\infty}^{-\alpha\Upsilon_0}- \int_{-\infty}^{-\alpha\Upsilon_0} q(r,x) dx.
\eqs
From the boundary conditions in \refe{lwr_bc}, we have
\bqs
xq(r,x)|_{-\infty}^{-\alpha\Upsilon_0}&=&0.
\eqs
Thus
\bqs
\pd{\Phi(k(r,x))}r&=&- \int_{-\infty}^{-\alpha\Upsilon_0} q(r,x) dx \leq 0,
\eqs
since $q(r,x)\geq 0$ in the LWR model. When $k(r,x)=k^*(x)$, $q(r,x)=0$, and $\pd{\Phi(k(r,x))}r=0$. In addition, if $k(r,x)\neq k^*(x)$ on a set of a positive measure, from \refe{lwr_kwt} $q(r,x)=\min\{d(r,x^-), s(r,x^+)\}>0$ for some $x$, and $\pd{\Phi(k(r,x))}r<0$.
\eop

\begin{theorem} If we use $\Phi(k(r,x))$ as the potential function,  the LWR model, \refe{lwr_sc}, is globally day-to-day asymptotically stable at the SPUE, according to  Lyapunov's second method.	
	\end{theorem}
{\em Proof}. From Theorem \ref{minimization_theorem}, we can see that the potential function reaches its minimum at $k^*(x)$. Further from Lemma \ref{derivative_lemma}, the potential function's derivative is strictly negative at  densities other than $k^*(x)$. Therefore, according to Lyapunov's second method \citep{LaSalle1960lyapunov},  the LWR model is asymptotically stable at the SPUE. \eop

Our study provides another interesting example to the literature of the global stability of hyperbolic conservation laws, which is rather scarce \citep{xu2002exponential,coron2007strict}.

\section{Conclusion}
In this study we first generalized the conceptual models for arrival time and scheduling payoff choices developed in \citep{jin2018_stable} for a single bottleneck with a generalized scheduling cost function. Then we presented the LWR model for the day-to-day traffic flow dynamics on an imaginary road for the scheduling payoff choice as well as the scheduling payoff user equilibrium (SPUE), which is equivalent to the departure time user equilibrium (DTUE) and the arrival time user equilibrium (ATUE). 
We further formulated a new optimization problem for the SPUE and demonstrated its equivalent to the optimization problem for the ATUE in \citep{iryo2007equivalent}. Finally we showed that the objective functions in the two optimization formulations can be used as the potential function for the LWR model and proved that the SPUE and, therefore, the DTUE and ATUE, are globally, asymptotically stable, by using Lyapunov's second method. From both the dynamical system and optimization formulations, we can see that the SPUE can be analytically solved, exist, and are unique.

In addition to offering a proof of the global stability of day-to-day departure time choice at a single bottleneck, this study also reveals the relationship between the dynamical system formulation and the optimization formulation. On the one hand, the objective function of the optimization formulation is the potential function of the dynamical system; on the other hand, the dynamical system model serves as a method to solve the optimization problem. Such a relationship has been observed in the corresponding dynamical system and optimization formulations of the static traffic assignment problem \citep[e.g.][]{jin2007dstap}, but is new for the departure time choice problem. 
The observation that the objective function of the optimization problems is the potential function for the dynamical system of day-to-day scheduling payoff choice suggests that the LWR model could be a behaviorally sound model to describe the day-to-day traffic flow of departure time choice on an imaginary road with increasing scheduling payoffs.

In the future we are interested in extending both dynamical system and optimization formulations for the SPUE in more complicated road networks, in which different vehicles can have different free-flow travel times, scheduling cost functions, values of time, and so on.
By studying the relationship between the two formulations,  we could obtain new insights regarding the choice behaviors and effective methods for solving the dynamic user equilibrium with simultaneous departure time and route choices. 
We will also be interested in calibrating and validating the behavioral principles incorporated into the model with real-world data. 
Such a globally stable behavioral model can provide more efficient departure time and route choice guidance for human drivers and connected and autonomous vehicles in more complicated networks.		
					
\section*{Acknowledgments}
Discussions with Prof. Kentaro Wada of the University of Tokyo are acknowledged. The views and results are the author's alone.

\pdfbookmark[1]{References}{references}

\end {document}